\documentclass[a4paper, 10pt]{article}

\usepackage{exscale}
\usepackage[centertags]{amsmath}
\usepackage{amssymb}
\usepackage{amsthm}
\usepackage{dsfont}
\usepackage{hyperref}
\usepackage[curve,matrix,arrow,ps]{xy}

\newtheorem{definition}{Definition}

\newtheorem{theorem}[definition]{Theorem}

\newtheorem{lemma}[definition]{Lemma}
\newtheorem{corollary}[definition]{Corollary}

\newcommand{\nd}{\noindent}

\newcommand{\dR}{{\mathbb R}}
\newcommand{\dC}{{\mathds C}}

\newcommand{\dN}{{\mathds N}}

\newcommand{\dZ}{{\mathds Z}}

\newcommand{\cC}{\mathcal{C}}
\newcommand{\cD}{\mathcal{D}}

\newcommand{\cF}{\mathcal{F}}
\newcommand{\cG}{\mathcal{G}}
\newcommand{\cH}{\mathcal{H}}
\newcommand{\cI}{\mathcal{I}}

\newcommand{\cK}{\mathcal{K}}
\newcommand{\cL}{\mathcal{L}}

\newcommand{\cN}{\mathcal{N}}
\newcommand{\cO}{\mathcal{O}}
\newcommand{\cP}{\mathcal{P}}

\begin{document}

\title{Rigid and Complete Intersection Lagrangian Singularities}
\author{Duco van Straten and Christian Sevenheck}
\date{November 27, 2003}
\maketitle

\begin{abstract}
In this article we prove a rigidity theorem for lagrangian singularities
by studying the local cohomology of the lagrangian de Rham complex that was introduced in \cite{SevStrat}. The  result can be applied to show the rigidity of all open
swallowtails of dimension $\ge 2$.
In the case of lagrangian complete intersection singularities
the lagrangian de Rham complex turns out to be perverse.
We also show that lagrangian complete intersection in dimension greater
than two cannot be regular in codimension one.
\end{abstract}
\renewcommand{\thefootnote}{}

\section{Introduction}

Since the work of Arnold and his school (\cite{ArnoldSwallowtail},
\cite{Arnold} and \cite{Giv3}), singular lagrangian subvarieties
in symplectic manifolds have become increasingly important in
different areas of mathematics. Arnold and Givental mainly studied
lagrangian projections and calculated normal forms for these
objects starting from the correspondence between such projections
from \emph{smooth} lagrangian germs to the base and generating
families. This does not, however, include the study of deformation spaces
which allows the lagrangian singularity itself to deform. In \cite{SevStrat}
 we considered the deformation problem
for a lagrangian singularity $(L,0)\subset (\dC^{2n},0)$ given by the
 deformation functor $\mathit{LagDef}^{loc}_{L,0}$
associating to a base space $S$ the set of isomorphism classes of
flat families $\cL \rightarrow S$ sitting inside $\dC^{2n}\times
S$ with the property that each fibre $\cL_s$ for $s\in S$ is
lagrangian in $\dC^{2n}\times \{s\}$. Similarly, one might define
a corresponding functor $\mathit{LagDef}_L$ for an analytic
lagrangian subspace $L$ inside a symplectic manifold $M$. The main
result of the quoted paper is a description of the tangent space
of this functor using the so-called lagrangian de Rham complex. We
recall this construction in section \ref{secLagdeRham} below.

In this paper we investigate some further properties of this
complex. We derive an inductive principle which can be used to
prove vanishing of the cohomology of the lagrangian de Rham
complex. This yields rigidity theorems for certain lagrangian
singularities of dimension higher than two and is similar in
spirit to the result of Schlessinger \cite{SchlessingerQuotient}
allowing to conclude that quotient singularities which are
regular in codimension two are rigid.
In \cite{SevStrat}, we also developed a constructive
method to calculate deformation spaces, but this was limited to
lagrangian surfaces. Therefore the results here are complementary
to our first paper, in that they extend the class of examples for
which deformations can be studied. On the other hand, the explicit
calculations from \cite{SevStrat} are used to make the induction
principle work.

The essential ingredients used in this article are the special
behavior of lagrangian deformations with respect to the canonical
stratification of a singularity and the local cohomology
of the lagrangian de Rham complex. One particular example of
lagrangian singularities to which our method applies are the
so-called open swallowtails. We show that they are all rigid.

The local cohomology sheaves of the lagrangian de Rham complex
also play a role in question of its perversity.  We show
here that lagrangian complete intersections have perverse
lagrangian de Rham complex. In this case, there is (via the
Riemann-Hilbert correspondence) a single $\cD$-module associated
to the lagrangian de Rham complex. This is consistent with an
abstract construction of this complex described in \cite{Thesis}.

A last result contained in this paper is concerned with the codimension
of the singular locus for lagrangian complete intersections.
We show that if such a singularity is regular in codimension two,
the tangent module is free. So the space is smooth for
all cases where the Zariski-Lipman problem is solved
affirmatively, in particular in the quasi-homogeneous case and the
case where the space is regular in codimension two.

\textbf{Acknowledgements}: We would like to thank A.~Givental for
calling our attention to his paper \cite{Giv5}.

\section{The lagrangian de Rham complex}
\label{secLagdeRham}

We recall in this section the construction from \cite{SevStrat} of
a sheaf complex associated to any Lagrangian variety.
The relationship of Lie algebroids and lagrangian singularities
is described in detail in \cite{Thesis}.
\begin{definition}
Let $L\subset \dC^{2n}$ be a lagrangian subvariety with defining
ideal sheaf $\cI\subset \cO_{\dC^{2n}}$. Denote by
$\cO_L:=\cO_{\dC^{2n}}/\cI$ the structure sheaf of $L$. The module
$\cI/\cI^2$ is the conormal module and has a structure of a Lie
algebroid over $\cO_L$, i.e., there are operations
$$
\{\,,\,\}: \cI/\cI^2 \times \cI/\cI^2 \longrightarrow \cI/\cI^2,\;\;\;\;\;\;\;\;\;\; \{\,,\,\}: \cI/\cI^2 \times \cO_L \longrightarrow \cO_L \\
$$
Define a sheaf complex $(\cC^\bullet_L, \delta)$, the lagrangian
de Rham complex by
$$
\cC^p_L:={\cH}\!om_{\cO_L}\left(\bigwedge^p\cI/\cI^2,\cO_L\right)
$$
and $\delta:\cC^p_L\rightarrow\cC^{p+1}_L$ with \small
$$
\begin{array}{l}
\left(\delta\left(\phi\right)\right)\left(h_1 \wedge \ldots \wedge h_{p+1}\right) := \\
  \sum_{i=1}^{p+1} \left(-1\right)^i
  \left\{h_i,\phi\left(h_1\wedge\ldots \wedge \widehat{h}_i \wedge \ldots
  h_{p+1}\right)\right\}
\\
 +  \sum\limits_{1 \leq i < j \leq p+1}\left(-1\right)^{i+j-1}
\phi\left(\left\{h_i,h_j\right\} \wedge  h_1  \wedge \ldots \wedge
\widehat{h}_i \wedge \ldots \wedge \widehat{h}_j \wedge \ldots
\wedge h_{p+1} \right)
\end{array}
$$
\normalsize
\end{definition}
We quote the main results from \cite{SevStrat} and \cite{Thesis}
concerning the lagrangian de Rham complex. The first one relates
$\cC^\bullet_L$ to the deformation theory of $L$.
\begin{theorem}
Consider the first three cohomology sheaves of $\cC^\bullet_L$.
Then
\begin{enumerate}
\item $\cH^0(\cC^\bullet_L)=\dC_L$ \item $\cH^1(\cC^\bullet_L)$ is
the sheaf of first order flat lagrangian deformations. This means
that at every point $p\in L$, the tangent space of the
functor $\mathit{LagDef}^{loc}_{L,p}$ is $H^1(\cC^\bullet_{L,p})$.
\item Let $(L,0)$ be either a complete
intersection or Cohen-Macaulay of codimension two. Suppose
moreover that $\cH^2(\cC^\bullet_L)=0$. Then the functor
$\mathit{LagDef}_{L,0}^{loc}$ is unobstructed.
\end{enumerate}
\end{theorem}
From the theory of Schlessinger, it is of obvious importance to
know whether the cohomology of the lagrangian de Rham complex is
finite. This is answered by the following result.
\begin{theorem}
\label{theoConstruct} Consider the canonical stratification of $L$
by embedding dimension, i.e., let $S^L_k:=\{p\in
L\,|\,\mathit{edim}_p(L)=2n-k\}$, where $k\in\{0,\ldots,n\}$.
Suppose that ``Condition P'' holds, that is, $\dim(S^L_k)\leq k$
for all $k$. Then the cohomology sheaves $\cH^p(\cC^\bullet_L)$
are constructible with respect to the canonical stratification. In
particular, for a germ $(L,0)$, $H^1(\cC^\bullet_{L,0})$ is a
finite dimensional vector space. Therefore, there is a formally
semi-universal deformation with respect to $\mathit{LagDef}_{L,0}^{loc}$.
\end{theorem}

\section{The rigidity theorem}

In this section, we state and prove our main theorem. The
technical tool used is the local cohomology of a sheaf, that is,
the derived functor of the functor $\Gamma_T(X,-)$ of sections of
a sheaf $\cF$ over a space $X$ with support in a closed subspace
$T$. Let us start with some preliminary lemmas.
In what follows we consider a lagrangian subvariety $X\subset \dC^{2n}$ which is
not necessarily Stein or contractible. $T\subset X$, $T \neq X$ is always a
closed analytic subspace.
\begin{lemma}
\label{lemDefLocCohHam} Denote by $\delta\cO_X\subset \cN_X$ the
image (sheaf) of the differential
$$
\delta:\cC^0_X=\cO_X\longrightarrow\cC^1_X=\cN_X
$$
Then we have
$$
H_T^0(\cH^1(\cC^\bullet_X))=\mathit{Ker}\left(H_T^1(\delta\cO_X)\rightarrow
H_T^1(\cN_X)\right)
$$
\end{lemma}
\begin{proof}
Consider the first three terms of the sheaf complex
$\cC^\bullet_X$ associated to the lagrangian subvariety
$X\subset\dC^{2n}$. It reads
$$
0\longrightarrow \cO_X\longrightarrow \cN_X \longrightarrow
\cC^2_X
$$
We know that
$\cH^0(\cC^\bullet_X)=\mathit{{\cK}er}(\cO_X\rightarrow\cN_X)=\dC_X$.
By splitting into short exact sequences, we obtain
$$
\begin{array}{c}
0 \longrightarrow \dC_X\longrightarrow \cO_X \longrightarrow \delta \cO_X \longrightarrow 0 \\ \\
0 \longrightarrow \delta \cO_X \longrightarrow \cK \longrightarrow \cH^1(\cC^\bullet_X) \longrightarrow 0 \\ \\
0 \longrightarrow \cK \longrightarrow \cN_X \longrightarrow \delta
\cN_X\longrightarrow 0
\end{array}
$$
Here $\cK=\mathit{{\cK}\!er}(\cN_X\rightarrow\cC^2_X)$ and
$\delta\cN_X=\mathit{{\cI}\!m}(\cN_X\rightarrow\cC^2_X)$. Now we
can apply the functor $H_T^\bullet(-)$ to each of these sequences.
This gives three long exact sequences of local cohomology sheaves.
However, we know in advance that sheaves of type
$\mathit{{\cH}\!om}_{\cO_X}(-,\cO_X)$ are torsion free, so in
particular $H_T^0(\cC^i_X)=0$ for all $i$. Moreover, $\dC_X$,
$\delta \cO_X$, $\cK$ and $\delta\cN_X$ are subsheaves of $\cO_X$,
$\cN_X$ resp. $\cC^2_X$, so for them the group $H^0_T(-)$ also
vanishes. We obtain exact sequences
$$
\begin{array}{c}
0 \longrightarrow H_T^1(\dC_X) \longrightarrow H_T^1(\cO_X)
\longrightarrow H_T^1(\delta \cO_X)
 \longrightarrow H_T^2(\dC_X) \\ \\
0 \longrightarrow H_T^0(\cH^1(\cC^\bullet_X)) \longrightarrow
H_T^1(\delta\cO_X)
 \longrightarrow H_T^1(\cK) \\ \\
0 \longrightarrow H_T^1(\cK) \longrightarrow H_T^1(\cN_X)
\longrightarrow H_T^1(\delta \cN_X)
\end{array}
$$
Combining the last two sequences yields the desired formula. The
first sequence will be used later.
\end{proof}
We need to investigate further the local cohomology of the sheaf
$\cH^1(\cC^\bullet_X)$.
\begin{lemma}
\label{lemLocSec} There is an exact sequence
$$
0 \longrightarrow H^0(X,\delta\cO_X) \longrightarrow
H^0(X\backslash T, \delta \cO_X) \longrightarrow
H^1_T(\delta\cO_X)
$$
If $X$ is Stein and contractible (e.g., a representative of a germ
$(X,0)$), then the last arrow in the above sequence is surjective.
\end{lemma}
\begin{proof}
Consider the following basic sequence in local cohomology (see
\cite{HartshorneGrothendieck}: Let $\cF$ be a sheaf on a
topological space $Y$ and $T$ any closed subspace, then:
\begin{equation}
\label{eqSequence} 0 \rightarrow H^0_T(\cF) \rightarrow H^0(Y,
\cF) \rightarrow H^0(Y\backslash T, \cF) \rightarrow H^1_T(\cF)
\rightarrow H^1(Y, \cF) \rightarrow \ldots
\end{equation}
For $Y=X\subset\dC^{2n}$ and $\cF=\delta\cO_X$, we know that
$H_T^0(\delta\cO_X)=0$. This gives the sequence in the general
case. Moreover, we can apply the usual cohomology functor to the
sequence
$$
0 \longrightarrow \dC_X\longrightarrow \cO_X \longrightarrow
\delta \cO_X \longrightarrow 0
$$
yielding
$$
\ldots \longrightarrow H^1(X, \cO_X) \longrightarrow H^1(X, \delta
\cO_X) \longrightarrow H^2(X, \dC_X) \longrightarrow \ldots
$$
In case that  $X$ is contractible ($H^2(X, \dC_X)=0$) and Stein
($H^1(X, \cO_X)=0$) the term $H^1(X,\delta\cO_X)$ vanishes.
\end{proof}
This last two results tell us how to understand sections of the
cohomology sheaf $\cH^1(\cC^\bullet_X)$ with support in a subspace
$T$, that is, deformations which do not deform the space
$X\backslash T$: these are elements of
$H^1_T(\delta\cO_X)$, thus, sections of $\delta\cO_X$
over $X\backslash T$ which do not extend over $T$. If we
consider the case $T=\mathit{Sing}(X)$, this means that a
deformation is trivial iff the hamiltonian vector field which
trivializes it on the regular part (because $\cH^1(\cC^\bullet_L)$
is zero on $X_{reg}$) extends over the whole of $X$.

\begin{theorem}
\label{theoMain} Let $L\subset \dC^{2n}$ be a representative of a
lagrangian singularity $(L,0)\subset (\dC^{2n},0)$ satisfying
Condition P. Denote by $S\subset L$ the singular locus. Let
$T\subset S$ be a closed analytic subspace in $L$ contained in the
singular locus. Suppose that
\begin{enumerate}
\item $H^1_T(\delta\cO_L)=0$
\item $H^0(L^*,\cH^1(\cC^\bullet_{L^*}))=0$, where $L^*:=L\backslash T$.
\end{enumerate}
Then $H^1(\cC_{L,0})=0$, i.e., $L$ is rigid under
lagrangian deformations.
\end{theorem}
\begin{proof}
Denote by $S^*$ the singular locus of $L^*$, obviously,
$S^*:=S\backslash T$. Note that $L_{reg}=L^*\backslash S^*$
because of $T\subset S$. From lemma \ref{lemLocSec}, applied to
the spaces $L$ and $L^*$, we obtain the following diagram
$$
\xymatrix@!0 { 0 \ar[rrr]  & & & H^0(L, \delta \cO_L)
\ar@{^{(}->}[dd]^\alpha \ar[rrrr] & & & & H^0(L_{reg}, \delta
\cO_L) \ar@{=}[dd] \ar[rrrr] && & & H^1_S(\delta
\cO_L) \ar[dd]^\beta \ar[rrr] & & & 0 \\ \\
0 \ar[rrr]  & & & H^0(L^*, \delta \cO_L) \ar[rrrr] & & & &
H^0(L_{reg}, \delta \cO_L) \ar[rrrr] && & & H^1_{S^*}(\delta
\cO_{L^*}) }
$$
Here $\alpha$ is the restriction map and $\beta$ is the induced map.
Moreover, a class $c\in
H^0(\cH^1(\cC^\bullet_L))=H^0_S(\cH^1(\cC^\bullet_L))$
corresponding to a flat lagrangian deformation of $L\subset
\dC^{2n}$ is represented by lemma \ref{lemDefLocCohHam} by a class
(denoted by the same letter) $c\in H^1_S(\delta\cO_L)$ which goes
to zero in $H^1_S(\cN_L)$. The same diagram, with the sheaf
$\delta\cO_L$ replaced by $\cN_L$ shows that $\beta(c)$ goes to
zero in $H^1_{S^*}(\cN_L)$. By lemma \ref{lemDefLocCohHam} we also know that
$$
H^0(\cH^1(\cC^\bullet_{L^*}))=H^0_{S^*}(\cH^1(\cC^\bullet_{L^*}))=
\mathit{Ker}\left(H^1_{S^*}(\delta\cO_L)\rightarrow
H^1_{S^*}(\cN_L)\right)
$$
which vanishes by the second hypothesis. So we get that
$\beta(c)=0$, this means that there is a section $\widetilde{c}$
extending $c$ over $L^*$.

We can apply lemma \ref{lemLocSec} again, this time to the pair
$(L,T)$, yielding the sequence
$$
0 \longrightarrow H^0(L,\delta\cO_L) \longrightarrow H^0(L^*,
\delta \cO_L) \longrightarrow H^1_T(\delta\cO_L)
$$
From the first hypothesis, we obtain that $\widetilde{c}$ extends
to the whole of $L$, which implies immediately that the original
class $c$ in $H^1_S(\delta\cO_L)$ is zero. Therefore, $L$ is
infinitesimal rigid.
\end{proof}
Using lemma \ref{lemDefLocCohHam}, the first condition implies in
particular that $H_T^0(\cH^1(\cC^\bullet_L))=0$, that is, there
are no deformations deforming only $T$. This is of course weaker
than the vanishing of $H^1_T(\delta\cO_L)$ but still sufficient:
By the same argument as above, we see that the class
$\widetilde{c}\in H^1_T(\delta\cO_L)$ maps to zero in
$H^1_T(\cN_L)$ thus defining an element in
$H_T^0(\cH^1(\cC^\bullet_L))$. But in applications, we will rather
prove that $H^1_T(\delta\cO_L)=0$, therefore, it is more
natural to impose this condition than the vanishing of
$H_T^0(\cH^1(\cC^\bullet_L))$.

In order to make use of these result, we have to find conditions
that give $H^1_T(\delta\cO_L)=0$ and $H^0(\cH^1(L^*, \cC^\bullet_{L*}))=0$.
We start with the first group. It sits in the exact sequence
$$
\ldots \longrightarrow H^1_T(\cO_L) \longrightarrow H^1_T(\delta
\cO_L) \longrightarrow H^2_T(\dC_L) \longrightarrow \ldots
$$
so a sufficient condition is the vanishing of the groups
$H^1_T(\cO_L)$ and $H^2_T(\dC_L)$. Obviously, $H^1_T(\cO_L)$ is of
analytic and $H^2_T(\dC_L)$ of topological nature.
\begin{lemma}
\label{lemVanishLocCoh} Let $\dim(L)\geq 2$ and $T$ be a closed
subspace such that $\textup{depth}(\cO_{L,0})\geq 2+\dim(T)$. Then
$H^1_T(\cO_L)=0$.
\end{lemma}
\begin{proof}
The well-known relation between local cohomology and $Ext$
leads to the statement that $H^p_T(\cF)=0$ is
equivalent to $\mathit{Ext}^p_{\cO_L}(\cG,\cF)=0$ for any sheaf
$G$ with $\mathit{supp}(\cG)\subset T$, see \cite{HartshorneGrothendieck},
proposition 3.7. By the lemma of Ischebeck (\cite{Matsumura}), $\mathit{Ext}^p_{\cO_L}(\cG,\cF)=0$
for all $p<\mathit{depth}(\cF)-\dim(\mathit{supp}(\cG))$. So for
$\cF=\cO_L$ and $\mathit{supp}(G)\subset T$ we obtain that
$H^1_T(\cO_L)=0$.
\end{proof}
The next step is to investigate the topological group
$H^2_T(\dC_L)$. First it follows from the sequence
\ref{eqSequence} that in case that $L$ is contractible (e.g, for a
representative of a germ $(L,0)$), we have
$H^2_T(\dC_L)=H^1(L\backslash T, \dC_L)$. The following lemma
lists some cases where the first homology of $L\backslash T$ is
zero.
\begin{lemma}
\label{lemVanishTopGroup} We consider a general situation of a
germ $(X,0)$ of a complex space.
\begin{enumerate}
\item Consider the normalization
$$
n:\widetilde{X} \longrightarrow X
$$
and suppose that $\widetilde{X}$ is smooth. Let $T$ a subspace of
codimension at least two such that $n$ induces an homeomorphism
from $\widetilde{X}\backslash \widetilde{T}$ to $X\backslash T$,
where $\widetilde{T} := n^{-1}(T)$. Then $H^1(X\backslash T, \dC)$
vanishes. \item Let $(X,0)$ be a rational normal surface
singularity and $T=\mathit{Sing}(X)=\{0\}$. Then we also have
$H^1(X\backslash T, \dC)=0$. \item Suppose that $X$ is a complete
intersection and $T$ a closed subspace of codimension at least
three which contains $\mathit{Sing}(X)$, then $H^1(X\backslash T,
\dC)=0$.
\end{enumerate}
\end{lemma}
\begin{proof}

\begin{enumerate}
\item This is obvious since $\widetilde{X}\backslash\widetilde{T}$
is simply connected and homeomorphic to $X\backslash T$. \item It
is known that the link $M$ of $(X,0)$ is a deformation retract of
$X\backslash T$. On the other hand, for rational singularities the
group $H^1(M,\dZ)$ is torsion (see, e.g.,
\cite{BrieskornRational}), so that $H^1(X\backslash T, \dC)$ is
zero. \item This can be found in \cite{GreuelThesis} or
\cite{Looijenga}. We sketch the argument: First it follows from a
result on the depth of the modules of differential forms on $X$
that
$$
H^1(\Omega^\bullet_{X,0},d)\cong H^1(\Gamma(X\backslash T,
\Omega^\bullet_{X\backslash T}),d)
$$
The same reasoning shows (using also the two spectral sequences
for the hypercohomology of a sheaf complex) that
$H^1(\Gamma(X\backslash T, \Omega^\bullet_{X\backslash T}),d)\cong
H^1(X\backslash T,\dC)$. By an analytic argument, one can show
that the de Rham complex of $X$ is exact in degree one. This
yields immediately that $H^1(X\backslash T,\dC)=0$.
\end{enumerate}
\end{proof}
Combining the last two lemmas, we get conditions for
$H^1_T(\delta\cO_L)$ to be zero. Whenever this is the case, a
lagrangian deformation of the germ $(L,0)$ comes (if it exist)
from a deformation of a transversal slice at a point $p\in
L \backslash T$. If we know that such deformations does not exist,
we can conclude that $L$ is rigid. This enables us for example
to conclude that any lagrangian rational triple point
in $\dC^4$ is rigid.
As a further consequence, we obtain from the third part of the last
lemma that lagrangian complete intersection singularities $L$ with
$\textup{codim}(\mathit{Sing}(L))>2$ are rigid. However, as we
will see in the last section, such objects simply do not exist.

\section{Applications}

We will use the theorem from the last section to prove rigidity
under lagrangian deformations of a number of examples including
the so-called open swallowtails. Givental introduces these
varieties in \cite{Giv5} as subvarieties of certain jet spaces in
order to obtain normal form results for systems of partial
differential equations. All examples studied in that paper are
obtained using \emph{generating functions} of special type. Recall
that for any function germ $F$ defined on a product of two smooth
spaces $B\times X$ such that the restriction $f$ of $F$ to
$\{0\}\times X$ defines a function germ with isolated critical
points, one can define (choosing coordinates $(x_1,\ldots,x_k)$ on
$X$ and $(q_1,\ldots,q_n)$ on $B$)
$$
Lag(F):=\{(p,q)\in T^*B\,|\,\exists x\in
X\;(\partial_{x_i}F)(x,q)=0\;;\;\;p_i=\partial_{q_i}F\;;\;\;\forall
i\} \subset T^*B
$$
It is well known that $Lag(F)$ is a lagrangian subvariety in
$T^*B$. Moreover, the generating function also gives rise to a
legendrian variety in $\dC^{2n+1}$ (with coordinates
$(u,\mathbf{p},\mathbf{q})$ and the standard contact structure
$u-\mathbf{p}\,d\mathbf{q}$), simply by setting $u=F(x,q)$. The
front of the lagrangian resp. legendrian variety is the image of
the projection to the $(u,\mathbf{q})$-space.

On the other hand, the space of polynomials
$$
\cP_{2n+1}= \left\{t^{2n+1}+\frac{a_1}{(2n-1)!}
t^{2n-1}+\frac{a_2}{(2n-2)!}t^{2n-2}+\ldots+a_{2n}\right\}
$$
carries a natural symplectic structure related to the
representation theory of $sl_2$. The subvariety consisting of all
polynomials having a root of multiplicity at least $n+1$ is
lagrangian and appears as generic singularity of the so-called
``obstacle problem'' (\cite{Giv3}). It is called $n$-dimensional
open swallowtail and was denoted $\Sigma_n$ in \cite{SevStrat}. We
will see that it can be described using generating functions. More
precisely, let $g_n(x,q):=x^{n+1}+q_1x^{n-1}+\ldots+q_n$ and set
$F_{n,k}(x,q):=\int_0^xg_n(s,q)^{k+1}ds$. Denote by $\Sigma_{n,k}$
the lagrangian subspace $Lag(F_{n,k})\subset\dC^{2n}$ and by
$\Lambda_{n,k}$ its front. The following lemma, extracted from
\cite{Giv3} and \cite{Giv5}, describes the geometry of the
singularities $\Sigma_{n,k}$ (and of its front $\Lambda_{n,k}$).
Some of these facts are needed later to apply our rigidity
theorem.
\begin{lemma}
\label{lemGeomLagVar}
\begin{enumerate}
\item Denote by $\cP_{m,n}$ the space of polynomials of degree
$(k+1)(n+1)+1$ with fixed highest coefficient, sum of roots equal
to zero and $n+1$ critical points of multiplicity $k+1$, i.e., all
polynomials of the form
$$
p_{q_1,\ldots,q_n,u}(s)=\int_0^x g_n(s,q_1,\ldots,q_n)^{k+1} ds -
u
$$
The front $\Lambda_{n,k}$ of the lagrangian singularity
$\Sigma_{n,k}$ is isomorphic to hypersurface of polynomials in
$\cP_{n,k}$ with multiple roots (such a root has automatically
multiplicity at least $k+2$). \item A smooth normalization of
$\Sigma_{n,k}$ is given by the map
$$
\begin{array}{rcl}
n\,:\,(\dC^n,0) & \longrightarrow & (\Sigma_{n,k},0) \\
(x,q_1,\ldots,q_{n-1}) & \longmapsto & (q_1,\ldots, q_n,
p_1,\ldots, p_n)
\end{array}
$$
here $q_n=x^{n+1}+\sum_{i=1}^{n-1}q_i x^{n-i}$,
$p_i:=\partial_{q_i}F_{n,k}$. \item The variety $\Sigma_{n,1}$ is
isomorphic to the $n$-dimensional open swallowtail $\Sigma_n$.
\item $(\Sigma_{n,k},0)$ is Cohen-Macaulay.
\end{enumerate}
\end{lemma}
\begin{proof}
\begin{enumerate}
\item This is almost a tautology: The front $\Lambda_{n,k}$ is the
graph of the generating function $F_{n,k}$, seen as a multi-valued
function (with $n+1$-sheets) on the base $B$. For any point
$\mathbf{q}=(q_1,\ldots,q_n)\in B$, let
$\lambda_1,\ldots,\lambda_{n+1}$ be the zeros of $g_n^{k+1}$. Then the
$n+1$ points of $\Lambda_{n,m}$ lying over $\mathbf{q}$ correspond
to the elements $p_{(\mathbf{q},u)} \in \cP_{n,k}$ with
$u=F(\lambda_i,\mathbf{q})$. Obviously, $\lambda_i$ is a zero of
$p_{(\mathbf{q},u)}$ and of its derivative, so belongs to the
discriminant in $\cP_{n,k}$. \item The map $n$ is generically one
to one and therefore the normalization. \item We will see that
$\cO_{\Sigma_n,0}$ and $\cO_{\Sigma_{n,1},0}$ can be identified as
subalgebras of their respective (smooth) normalization. Following
\cite{Giv3}, the normalization of $\Sigma_n$ is given by the
following map
$$
\begin{array}{rcl}
\varphi:\widetilde{\Sigma}_n\cong\dC^n & \longrightarrow & \Sigma_n \subset \cP_{2n+1} \\
(x, a_1, \ldots, a_{n-1}) & \longmapsto &
(t-x)^{n+1}\cdot(t^n+b_1t^{n-1}+\ldots+b_{n-1})
\end{array}
$$
where $b_i\in\cO_{\widetilde{\Sigma}_n,0}$ are chosen such that
the coefficient of $t^{2n+1-i}$ in the polynomial
$\varphi(x,\mathbf{a})$ is precisely $a_i/(2n+1-i)!$ for
$i=1,\ldots,n-1$ (in particular, $b_1=(n+1)t$). Then we get
$$
\cO_{\Sigma_n,0} =
\left\{f\in\dC\{t,a_1,\ldots,a_{n-1}\}\;|\;f=\int_0^x\!\!
Q(s,\mathbf{a})F_n(s,\mathbf{a})ds+C(\mathbf{a})\right\}
$$
On the other hand, it is shown in \cite{Giv5} that
$$
\cO_{\Sigma_{n,k},0}=\left\{f\in\dC\{x,q_1,\ldots,q_{n-1}\}\;|\;f=\int_0^x\Phi(s,\mathbf{q})g_n(s,q)^k
ds +Q(\mathbf{q})\right\}
$$
So $\cO_{\Sigma_{n,\!1},\,0}\cong\cO_{\Sigma_n,0}$. \item One has
to show that the finite analytic mapping $
(\Sigma_{n,k},0)\rightarrow(B,0)$ makes $\Sigma_{n,k}$ into a
\textbf{free} $\cO_{B,0}$-module of rank $n+1$. This is done in
\cite{Giv5} (for $k=1$, this map is simply $n$-fold
differentiation). Then the statement follows.
\end{enumerate}
\end{proof}

From the first point of the lemma, we deduce
\begin{lemma}
Let
$$
\{0\}\subset\Sigma_{n,k}^{(1)}\subset\ldots\subset\Sigma_{n,k}^{(n-1)}
\subset \Sigma_{n,k}^{(n)}=\Sigma_{n,k}
$$
be the canonical stratification with $\dim(\Sigma_{n,k}^{(k)})=k$
(Condition P). Let $p\in \Sigma_{n,k}^{(i)}\backslash
\Sigma_{n,k}^{(i-1)}$, then we have $(\Sigma_{n,k}, p) \cong
(\Sigma_{n-i,k},0) \times(\dC^i,0)$.
\end{lemma}
\begin{proof}
That $\Sigma_{k,n}$ locally decomposes into a product of a
lagrangian variety and a smooth germ is a general fact (this is
the essential ingredient in the proof of theorem
\ref{theoConstruct}, see \cite{SevStrat} and \cite{Thesis}). We
only need to show that the transversal section is precisely
$(\Sigma_{n-i,k},0)$. First it is obviously sufficient to do case
$i=n-1$. For this, we will show that the transversal singularity
of the front $\Lambda_{n,k}$ is $\Lambda_{n-1,k}$. This follows
directly from the description of the front given as discriminant
in the polynomial space $\cP_{n,k}$. A general polynomial $P$ in
this space can be written in the form
$$
\int_0^x(s-\lambda)^{k+1}(s-\mu)^{k+1}(s^{n-1}+(\lambda+\mu)
s^{n-2}+q'_1s^{n-3}\ldots+q'_{n-2})^{k+1}ds
$$
with the additional condition that there is a common zero of $P$
and its derivative. If $\lambda=\mu$, then the polynomial $P$
represents a point $\widetilde{p} \in \Lambda_{n,k}$ corresponding
to the point $p\in\Sigma_{n,k}$ from above. A transversal section
at $\widetilde{p}$ is given (in appropriate local coordinates) by
setting $\lambda=\textup{const}$ and by translating the argument.
Therefore, in a neighborhood of $\widetilde{p}$ a point of such a
transversal section is represented as
$\int_0^x(s-\mu)^{k+1}(s^{n-1}+\mu s^{n-2}+
\widetilde{q}_1s^{n-3}\ldots+\widetilde{q}_{n-2})^{k+1}ds$, that
is, corresponds to a point in $\Lambda_{n-1,k} \subset
\cP_{n-1,k}$.
\end{proof}
In \cite{SevStrat}, an algorithm to calculate
$H^1(\cC^\bullet_{L,0})$ for quasi-homogenous lagrangian surface
singularities was described. For the spaces $\Sigma_{2,k}$ one
obtains by computer calculation
\begin{lemma}
\label{lemDirectCalc} $H^1(\cC^\bullet_{\Sigma_{2,k},0})=0$ for
$k=2,3,4,5$. In these cases, as in the examples studied in
\cite{SevStrat} the spectral numbers of the local system
$\cH^1(\cC^\bullet_L)_{|Sing(L)}$ (for $L=\Sigma_{2,k}$) have a
symmetry property.
\end{lemma}
For higher $k$ the computation is possible in the same way and
limited only by computer power. Conjecturally, all $\Sigma_{2,k}$
are rigid. For all $k$ such that $\Sigma_{2,k}$ is rigid, we can use theorem
\ref{theoMain} to obtain.
\begin{theorem}
Suppose that for fixed $k$, the lagrangian singularity
$(\Sigma_{2,k},0)\subset (\dC^4,0)$ is rigid. Then for all $n>2$
$(\Sigma_{n,k},0)\subset(\dC^{2n},0)$ is rigid.
\end{theorem}
\begin{corollary}
All open swallowtails of dimension greater one are rigid
lagrangian singularities.
\end{corollary}
\begin{proof}[Proof of the theorem]
We do induction on $n$. For $n=2$, we are done by hypothesis.
Otherwise, we know that for
$p\in\Sigma_{n,k}^{(1)}\backslash\{0\}$, there is a decomposition
$(\Sigma_{n,k},p)\cong(\Sigma_{n-1,k},0)\times(\dC,0)$ and moreover,
$H^1(\cC^\bullet_{\Sigma_{n,k},p})\cong
H^1(\cC^\bullet_{\Sigma_{n-1,k},0})$. This last group is zero by the
induction hypothesis. This implies that for $T=\{0\}\subset
\Sigma_{n,k}$, we have $H^0(\Sigma_{n,k}\backslash T,
\cH^1(\cC^\bullet_{\Sigma_{n,k}\backslash T}))=0$. The second
point we need to check in order to apply theorem \ref{theoMain} is
the vanishing of $H^1_T(\delta\cO_{\Sigma_{n,k}})$. We use lemma
\ref{lemVanishLocCoh}: We need that $T$ is of codimension at least
two and that $\textup{depth}(\cO_L)>\dim(T)+2$ which is obviously
satisfied in view of the last point of lemma \ref{lemGeomLagVar}.
Moreover, the second statement of this lemma gives smoothness of
the normalization of $(\Sigma_{n,k},0)$, so that the second
(topological) condition of lemma \ref{lemVanishLocCoh} is also
satisfied. Therefore, $H^1_T(\delta\cO_{\Sigma_{n,k}})=0$. Now we
can apply theorem \ref{theoMain}, which proves rigidity of
$\Sigma_{n,k}$.
\end{proof}

\section{Lagrangian complete intersections}

The perversity condition mentioned in the introduction involves
study of the local cohomology of the lagrangian de Rham complex.
For that reason, it is quite natural to include it here. It turns
out that a positive answer to this problem is possible in the case
of lagrangian complete intersections. Let us first recall what it
means for a complex to be perverse. Consider a, say, complex space
$X$ of dimension $n$ and a sheaf complex $\cK^\bullet$ on $X$ (we
suppose for simplicity that it is concentrated in non-negative
degrees). Then there are two condition, called first and second
perversity conditions. The first one states that
$$
\dim supp(\cH^i(\cK^\bullet)) \leq n-i
$$
for all $i\leq 0$. The second one (also called co-support condition) involves the derived functor
$\dR\Gamma_T$ (seen as functor in the derived category), where $T$
is a closed analytic subspace in $X$. It states that
$$
\dim supp(\dR^q\Gamma_T(\cK^\bullet)) < \dim(T)
$$
for any such $T$ and for all $i\in\{0,\ldots,n-\dim(T)-1\}$. We
also recall the spectral sequence with $E_2$-term
$\cH^p(\cH_T^q(\cK^\bullet))$ and which converges to
$\dR^{p+q}\Gamma_T(\cK^\bullet)$. Now consider the case
$\cK^\bullet=\cC^\bullet$.
\begin{theorem} \label{theoCoSupport}
Let $L\subset \dC^{2n}$ be a representative of a lagrangian
complete intersection singularity. Then the complex
$\cC^\bullet_L$ is perverse.
\end{theorem}
\begin{proof}
The first condition is easily verified using the decomposition of
a lagrangian variety around a point of non-maximal embedding
dimension (this has already been done in \cite{SevStrat}).

Consider the above spectral sequence. $L$ is a complete
intersection, therefore, the conormal module and hence the modules
$\cC^p_L$ are locally free. In particular,
$\textup{depth}(\cC^p_L)=n$. By the lemma of Ischebeck (see the
proof of lemma \ref{lemVanishLocCoh}), we have that
$\cH_T^q(\cC^p_L)=0$ for all $q<n-\dim(T)$. This implies the
vanishing of the corresponding local hypercohomology
$\dR^q\Gamma_T(\cC^\bullet_L)$, as required.
\end{proof}
\begin{corollary}
Let $L$ be lagrangian with $\dim(L)\leq 3$ and
$\textup{depth}(L)\geq 2$. Then $\cC^\bullet_L$ is perverse.
\end{corollary}
\begin{proof}
The proof of the last theorem shows that whenever we have a
vanishing of $\cH^q(\cC^p_L)$, we get vanishing of the
hypercohomology. But there is a general statement (see, e.g.,
\cite{SchlessingerQuotient}) that for a space $X$ of depth at
least two, sheaves of type $\mathit{{\cH}\!om}_{\cO_X}(\cF,\cO_X)$
are of depth at least two. Up to dimension three, this vanishing
result is sufficient for the co-support condition to be satisfied.
\end{proof}
The natural question whether there exist examples of lagrangian
singularities with non-perverse lagrangian de Rham complex is
still open. If one looks at the open swallowtail
$\Sigma_4\subset\cP_9\cong\dC^8$, it would be sufficient to have
$\textup{depth}(\cN_{\Sigma_2})=2$ in order to get a counterexample, but we were not able to compute this depth.

We remark that the co-support condition simplifies due to the
decomposition principle as follows: Let $\{0\} = L^{(0)} \subset
L^{(1)}\subset\ldots\subset L^{(n)}=L$ be the canonical
stratification. Then it is sufficient to show the co-support
condition only for subspaces $T=L^{(i)}$. Moreover, if we can show
that for all $i$ and all $q\in\{0,\ldots,n-i-1\}$,
$$
supp(\dR^q\Gamma_{L^{(i)}}(\cC^\bullet_L)) \subset L^{(i-1)}
$$
(where $L^{(-1)}:=\emptyset$), then we are done by ``Condition
P''. This amounts to show that for $p\in L^{(i)} \backslash
L^{(i-1)}$, the stalk $\dR^q\Gamma_{L^{(i)}}(\cC^\bullet_L)_p$ is
zero. But we know that $(L,p)\cong(L',p')\times(\dC^i,0)$ with
$p'\in L^{(0)}$ and that $\cC^\bullet_{L,p}$ is quasi-isomorphic
to $\pi^{-1}\cC^\bullet_{L',p'}$ (with
$\pi:(L,p)\rightarrow(L',p')$ the projection). Therefore
$$
\dR^q\Gamma_{L^{(i)}}(\cC^\bullet_L)_p =
\dR^q\Gamma_{L^{(i)}}(\pi^{-1}\cC^\bullet_{L'})_p=
\dR^q\Gamma_{{L'}^{(0)}}(\cC^\bullet_{L'})_{p'}
$$
So if we know for a class of lagrangian singularities that the
transversal slices also belongs to this class (as, e.g., for
complete intersections), it suffices to show that
$\dR^q\Gamma_{\{0\}}(\cC^\bullet_L)=0$ for all $0\leq q <n$ for
all $L$ in this class.

We add here a statement giving a partial answer to a question
on the singular locus of lagrangian complete intersections.
\begin{theorem}
Let $(L,0)\subset(\dC^{2n},0)$ be a lagrangian complete
intersection singularity such that
$\textup{codim}(\mathit{Sing}(L)) \ge 2$.
Then the tangent module $\Theta_{L,0}$ is free.
\end{theorem}

\begin{proof}[Proof of the theorem]
Let $I\subset\cO_{\dC^{2n},0}$ be the defining ideal of $(L,0)$.
From \cite{SevStrat}, we have the following diagram
$$
\xymatrix@C=0.5cm {
&{I/I^2} \ar[r] \ar[d]^{\alpha} & \Omega^1_{\dC^{2n},0} \otimes \cO_{L,0} \ar[d]^{\cong} \ar[r] & \Omega^1_{L,0} \ar[d]^{\widetilde{\alpha}} \ar[r] & 0\\
0 \ar[r] &\Theta_{L,0} \ar[r] & \Theta_{\dC^{2n},0} \otimes
\cO_{L,0} \ar[r] & N_{L,0} \ar[r] & T^1_{L,0} \ar[r] &0 }
$$
where $\alpha$ and $\widetilde{\alpha}$ are isomorphisms on
$L_{reg}$. By the snake lemma,
$\mathit{Coker}(\alpha)\cong\mathit{Ker}(\widetilde{\alpha})$. We
know that $N_{L,0}$ is torsion free and that the kernel of
$\widetilde{\alpha}$ is concentrated on the non-smooth locus,
hence, $\mathit{Coker}(\alpha)\cong
\mathit{Tors}(\Omega^1_{L,0})$. Now if $L$ is a complete
intersection then it follows from \cite{GreuelThesis} that
$\Omega_{L,0}^p$ is torsion free for all
$p<\textup{codim}(\mathit{Sing}(L))$, in particular,
$\Omega_{L,0}^1$ is torsion free under the hypotheses of the
theorem. This shows that $\alpha$ is surjective. For a complete
intersection, $I/I^2$ is free and $I/I^2 \rightarrow
\Omega^1_{\dC^{2n},0}\otimes \cO_{L,0}$ is injective. Therefore,
$\Theta_{L,0}\cong I/I^2$ is free.
\end{proof}

From the freeness of $\Theta_{L,0}$ one would like to conclude
that $(L,0)$ is in fact smooth. This is the celebrated
\emph{Zariski-Lipman}-conjecture.
\begin{center}
\textit{Let $R$ be an analytic $\dC$-algebra such that the
$R$-module \\
$\Theta_R:=\mathit{Der}_{\dC}(R,R)$ is free. Then $R$ is smooth.}
\end{center}
This conjecture is proved in a number of cases. The first case is the graded one, due to Platte \cite{Platte}, starting from a proof in the algebraic
case by Hochster, \cite{Hochster1}, \cite{Hochster2}.
\begin{lemma}
Let $A$ be a positively graded analytic algebra, that is, there is
$E\in\Theta_A$ such that the maximal ideal $\mathbf{m}_A$ is
generated by elements $x_i$ with $E(x_i)=w(i)$, where $w(i) \in
\dN_{>0}$. If $\Theta_A$ is a free $A$-module, then $A$ is
regular.
\end{lemma}
If $R$ is not graded, one has to use rather different techniques.
The following lemma (\cite{vanStratenSteenbrink}) relates the
Zariski-Lipman conjecture with the question of extendability of
differential forms on $R$ to its resolution.
\begin{lemma}
Let $(X,0)$ the germ of an analytic space $X$. Consider a
resolution $\pi:\widetilde{X}\rightarrow X$ with
$\pi_*\Theta_{\widetilde{X}} \cong \Theta_X$. Let
$U:=X\backslash\mathit{Sing}(X)$. If the natural morphism
$\Omega_{\widetilde{X}} \rightarrow \pi^*\Omega_U$ is surjective,
then $X$ is smooth if $\Theta_X$ is locally free.
\end{lemma}
\begin{proof}
The idea is simply that a basis $\theta_1,\ldots,\theta_n$ of
$\Theta_X$ gives rise to vector fields on
$\widetilde{\theta}_1,\ldots,\widetilde{\theta}_n$ on
$\widetilde{X}$ tangent to the exceptional locus $E$ of the
resolution. On $\widetilde{X}\backslash E$, there are independent
forms $\alpha_1,\ldots,\alpha_n$ dual to these vector fields which
extends over $E$. This is a contradiction, as for any point $p\in
E$, the vectors $\widetilde{\theta}_i(p)$ cannot be linearly
independent, because $\dim(E) < n$, unless $X$ is smooth.
\end{proof}
In the quoted paper, the extendability of differential $p$-forms
on isolated singularities is studied and the authors prove that
any $p$-form with $p<\dim(R)-1$ is extendible. Flenner
(\cite{Flenner}) showed that more generally, for any space $X$, a $p$-form
on $X\backslash \mathit{Sing}(X)$ with $p < \textup{codim}(\mathit{Sing}(X))-1$
extends to a resolution $\widetilde{X}$ of $X$. Therefore, one has
\begin{corollary}
Let $R$ be any analytic algebra such that $\textup{codim}(\mathit{Sing}(R))\geq 3$.
Then the Zariski-Lipman conjecture is true.
\end{corollary}

One of the sources of lagrangian singularities are Frobenius manifolds, where
they arise as \emph{spectral covers} of the multiplication on the tangent
bundle. It was asked in \cite{Hertling}, chapter 14, if there exist
an isolated Gorenstein, hence complete intersection, lagrangian
surface singularity. In the quasi-homogeneous case this is excluded by
our theorem.  The case of a non-quasi-homogenous lagrangian isolated
complete intersection surface singularity remains open, because
the Zariski-Lipman conjecture is unproven in this key case.

\bibliographystyle{amsalpha}
\providecommand{\bysame}{\leavevmode\hbox to3em{\hrulefill}\thinspace}
\providecommand{\MR}{\relax\ifhmode\unskip\space\fi MR }
\providecommand{\MRhref}[2]{%
  \href{http://www.ams.org/mathscinet-getitem?mr=#1}{#2}
}
\providecommand{\href}[2]{#2}

\vspace*{1cm}

\nd
Christian Sevenheck\hfill Duco van Straten\\
Ecole Normale Supérieur \hfill  Johannes-Gutenberg-Universit\"at Mainz \\
Département de mathématiques et applications \hfill FB 17, Mathematik\\
45, rue d'Ulm \hfill 55099 Mainz\\
75230 Paris cedex 05, France \hfill Germany\\
Christian.Sevenheck@ens.fr \hfill
straten@mathematik.uni-mainz.de

\end{document}